\documentclass[12pt, reqno]{amsart}
\usepackage{amsmath, amsthm, amscd, amsfonts, amssymb, graphicx, color}
\usepackage[bookmarksnumbered, colorlinks, plainpages]{hyperref}

\newtheorem{theorem}{Theorem}[section]

\newtheorem{corollary}[theorem]{Corollary}
\theoremstyle{definition}

\theoremstyle{remark}

\numberwithin{equation}{section}
\begin{document}
\begin{center}
{\large{\textbf{Darboux Rectifying curves on a smooth surface}}}
\end{center}
\vspace{0.1 cm}
\begin{center}
{Buddhadev Pal$\footnote{Corresponding author}$, Akhilesh Yadav }
\end{center}
\vskip 0.3cm
\begin{center}
Department of Mathematics\\
Institute of Science\\
Banaras Hindu University\\
Varanasi-221005, India.\\

\vskip 0.3cm
E-mail: pal.buddha@gmail.com\\
E-mail: akhilesha68@gmail.com\\
\end{center}
\vskip 0.5cm
\textbf{Abstract:}
The main aim of this paper is to investigate Darboux rectifying curves on a smooth surface immersed in the Euclidean space. First, we discuss the component of the position vector of a Darboux rectifying curve on a smooth immersed surface under the isometry of surfaces. Next we find a sufficient condition for the conformal invariance of Darboux rectifying curve.
\section{\textbf{Introduction}}
In 2003, B. Y. Chen \cite{che} provided the concept of rectifying curve in $\mathbb{R}^3$. The space curves whose position vector always lies in its rectifying plane, called rectifying curve. To know more about the characterization of rectifying curve we refer the reader to see (\cite{cd}, \cite{dca}, \cite{ygy}). Recently A. A. Shaikh and P. R. Ghosh \cite{sg} discussed about rectifying curves on a smooth surface. Also from the books A. Pressley \cite{ap} and M. P. do Carmo \cite{mp} gave the details idea about characterization of space curve. When we study geometry of space curves on a smooth surface, arrived Frenet frame and as well as arrived Darboux frame, exists at any non-umbilic point of a surface embedded in Euclidean space. It is named after French mathematician Jean Gaston Darboux. In (\cite{uk}, \cite{hh}) authors have studied Frenet and Darboux rotation vectors of curves on Time-like surfaces. They have given the Darboux frame of the curves according to the Lorentzian characters of surfaces and the curves.\\\\
We know that for the case of isometry between two surfaces, length of the curves and angle between the intersecting curves  are preserved and in the case of conformal motion, only angles are preserved and not necessarily distances. Now from \cite{hgk}, we can say isometric maps are a subset of conformal maps with the dilation factor is equal to 1.\\\\
Motivated from the above papers we study Darboux rectifying curve on a smooth surface and find sufficient condition for invariance of Darboux rectifying curve under isometry and under conformal map. We organize our paper as follows: In section 2, we discuss the relations between Frenet and Darboux frames in Euclidean space. In section 3, we investigate Darboux rectifying curve on a smooth surface immersed in the Euclidean space. Here we prove that under isometry of surfaces the components of the position vector of a Darboux rectifying curve on a smooth surface along the normal to the curve $(\kappa(s) = \bar{\kappa}(s))$ and along any tangent vector on the surface are invariants. In the last section, we study the conformal properties of a Darboux rectifying curve which totally depends upon its position vector field.
\section{\textbf{Preliminaries}}
Let $\gamma = \gamma(s)$ be a unit speed curve in $\mathbb{R}^3$, parameterized by its arc length $s$. Let $T$, $N$, $B$ be the tangent, principal normal and binormal vector fields along the curve $\gamma(s)$. Then the Frenet formulae is given by
\[\left( \begin{array}{c}
          T^{'} \\
          N^{'}\\
          B^{'} \end{array} \right) = \left( \begin{array}{ccc}
           0 & \kappa & 0 \\
           -\kappa & 0 & \tau\\
           0 & -\tau & 0 \end{array} \right)\left( \begin{array}{c}
           T \\
           N\\
           B \end{array} \right), \]
where the functions $\kappa(s)>0$ and $\tau(s)$ are called the curvature and the torsion of the curve and a prime denotes differentiation with respect to $s$.\\\\
Let $\{T, U, P\}$, Darboux frame of $\gamma(s)$, where $T$ as the tangent vector of $\gamma$ and $U$ is the unit normal to the surface $M$ and $P = U\times T$, then the relation between Frenet and Darboux frames takes the form (\cite{mp}, \cite{bo}):
\[\left( \begin{array}{c}
          T \\
          P\\
          U \end{array} \right) = \left( \begin{array}{ccc}
           1 & 0 & 0 \\
           0 & \cos\alpha & \sin\alpha\\
           0 & -\sin\alpha & \cos\alpha \end{array} \right)\left( \begin{array}{c}
           T \\
           N\\
           B \end{array} \right), \]
where there is a unique angle $\alpha$ such that a rotation in the plane $N$ and $B$ produces the pair $P$ and $U$. Therefore the Darboux formulae of $\gamma(s)$ is given as follows:
\[\left( \begin{array}{c}
          T^{'} \\
          P^{'}\\
          U^{'} \end{array} \right) = \left( \begin{array}{ccc}
           0 & \kappa_{g} & \kappa_{n} \\
           -\kappa_{g} & 0 & \tau_{g}\\
           -\kappa_{n} & -\tau_{g} & 0 \end{array} \right)\left( \begin{array}{c}
           T \\
           P\\
           U \end{array} \right), \]
where $\kappa_{g}$, $\kappa_{n}$ and $\tau_{g}$ are the geodesic curvature, normal curvature and geodesic torsion, respectively.\\\\
Let $\gamma: (\alpha, \beta) \rightarrow M$ be a unit speed parameterized curve on the coordinate chart $\eta: V \rightarrow M$  of the smooth surface $M$. Thus the curve $\gamma(s)$ is contained in the image of a surface patch $\eta$. Hence we can write
\begin{equation}
\gamma(s) = \eta(u(s), v(s)).\end{equation}
Now, differentiating equation (2.1) with respect to $s$, we get
\begin{equation}
T(s) = \gamma'(s) = \eta_{u}u' + \eta_{v}v',
\end{equation}
\begin{equation}
 \Rightarrow T'(s)= \eta_{u}u'' + \eta_{v}v'' + u'^{2}\eta_{uu} + 2u'v'\eta_{uv} + v'^{2}\eta_{vv}.
\end{equation}
If $\kappa(s)$ is the curvature of $\gamma(s)$ and $U$ is the unit normal to $M$ then the principal normal vector $N(s)$ of the curve at the point $\gamma(s)$ is given by
\begin{equation}
 N(s)= \frac{1}{\kappa(s)}(\eta_{u}u'' + \eta_{v}v'' + u'^{2}\eta_{uu} + 2u'v'\eta_{uv} + v'^{2}\eta_{vv}),
\end{equation}
and binormal vector $B(s)$ is given by
 $$B(s)= \frac{1}{\kappa(s)}[(\eta_{u}u' + \eta_{v}v')\times (\eta_{u}u'' + \eta_{v}v'' + u'^{2}\eta_{uu} + 2u'v'\eta_{uv} + v'^{2}\eta_{vv})],$$
\begin{eqnarray}
i.e., B(s)= \frac{1}{\kappa(s)}[(u'v''-v'u'')U + u'^{3}\eta_{u}\times \eta_{uu} + 2u'^{2}v'\eta_{u}\times \eta_{uv} + u'v'^{2}\eta_{u}\times \eta_{vv}\nonumber\\ + u'^{2}v'\eta_{v}\times \eta_{uu} + 2u'v'^{2}\eta_{v}\times \eta_{uv} + v'^{3}\eta_{v}\times \eta_{vv}].
\end{eqnarray}
Also $U$ is the unit normal to the surface, therefore
\begin{equation}
 U(s)= \frac{\eta_{u}\times \eta_{v}}{||\eta_{u}\times \eta_{v}||} = \frac{\eta_{u}\times \eta_{v}}{\sqrt{EG - F^{2}}},
\end{equation}
and
$$ P(s)= U(s)\times T(s) = \frac{\eta_{u}\times \eta_{v}}{||\eta_{u}\times \eta_{v}||}\times (\eta_{u}u' + \eta_{v}v'),$$
\begin{equation}
 \Rightarrow P(s) = \frac{1}{||\eta_{u}\times \eta_{v}||}(Eu'\eta_{v} + F(v'\eta_{v} - u'\eta_{u}) - Gv'\eta_{u}).
\end{equation}
where $E =  \eta_{u}.\eta_{u}$, $F = \eta_{u}.\eta_{v}$ and $G = \eta_{v}.\eta_{v}$ are the coefficients of the first fundamental form.
\section{\textbf{Darboux rectifying curves on a smooth surface}}
In this section, we study Darboux rectifying curve on a smooth surface and find sufficient condition for invariance of Darboux rectifying curve under isometry.\\\\
\textbf{Definition 3.1:} {\em A curve $\gamma(s)$ on a smooth surface $M$ whose position vector lies in the \{T,P\}-Darboux rectifying plane is called Darboux rectifying curve}.\\\\
Thus the equation of Darboux rectifying curve is given by:
\begin{equation}
\gamma(s) = \lambda(s)T(s) + \mu(s)P(s),
\end{equation}
for smooth functions $\lambda(s)$ and $\mu(s)$.\\\\
Now by using equations (2.2) and (2.7) in (3.1), we get
\begin{equation}
\gamma(s) = \lambda(s)(\eta_{u}u' + \eta_{v}v') + \frac{\mu(s)}{\sqrt{EG - F^{2}}}(Eu'\eta_{v} + F(v'\eta_{v} - u'\eta_{u}) - Gv'\eta_{u}).
\end{equation}
Next we consider the expression of the derivative map of $\gamma(s)$ i.e.,$f_{*}(\gamma(s))$ as product of a 3 $\times$ 3 matrix $f_{}*$ and a 3 $\times$ 1 matrix $\gamma(s)$.
\begin{theorem}
Let $f$ be an isometry between two smooth surfaces $M$ and $\bar{M}$. Suppose $\gamma(s)$ be a Darboux rectifying curve on $M$. Then $\bar{\gamma}(s)$ is a Darboux rectifying curve on $\bar{M}$ if $$\bar{\gamma}(s) = f_{*}(\gamma(s)).$$
\end{theorem}
\begin{proof}
Let $f$ be an isometry between two smooth surfaces $M$ and $\bar{M}$. Then, $f_{*}: T_{p}M \rightarrow T_{f(p)}\bar{M}$ such that $f_{*}\eta_{u} = \bar{\eta_{u}}$ and $f_{*}\eta_{v} = \bar{\eta_{v}}$. Also since $M$ and $\bar{M}$ are isometric, $E = \bar{E}$, $F = \bar{F}$ and $G = \bar{G}$.\\\\ Suppose $\bar{\gamma(s)} = f_{*}(\gamma(s))$.\\\\
Now from the equation (3.2) of Darboux rectifying curve, we have
\begin{align*}
\bar{\gamma(s)} = \lambda(s)(u'f_{*}(\eta_{u}) + v'f_{*}(\eta_{v})) + \frac{\mu(s)}{||f_{*}\eta_{u}\times f_{*}\eta_{v}||}\\(Eu'f_{*}\eta_{v} + F(v'f_{*}\eta_{v} - u'f_{*}\eta_{u}) - Gv'f_{*}\eta_{u}),
\end{align*}
$$\Rightarrow \bar{\gamma(s)} = \lambda(s)(u'\bar{\eta_{u}} + v'\bar{\eta_{v}}) + \frac{\mu(s)}{||\bar{\eta_{u}}\times \bar{\eta_{v}}||}(Eu'\bar{\eta_{v}} + F(v'\bar{\eta_{v}} - u'\bar{\eta_{u}}) - Gv'\bar{\eta_{u}}),$$
$$\Rightarrow \bar{\gamma(s)} = \lambda(s)(u'\bar{\eta_{u}} + v'\bar{\eta_{v}}) + \frac{\mu(s)}{\sqrt{\bar{E}\bar{G} - \bar{F}^{2}}}(\bar{E}u'\bar{\eta_{v}} + \bar{F}(v'\bar{\eta_{v}} - u'\bar{\eta_{u}}) - \bar{G}v'\bar{\eta_{u}}),$$
$$\Rightarrow \bar{\gamma(s)} = \bar{\lambda}(s)\bar{T}(s) + \bar{\mu}(s)\bar{P}(s),$$\\
for smooth functions $\bar{\lambda}(s)$ and $\bar{\mu}(s)$. Thus $\bar{\gamma}(s)$ is a Darboux rectifying curve on $\bar{M}$.
\end{proof}
\textbf{Note:} From the Theorem 3.1, we see that the functions $\lambda(s)$, $\mu(s)$  and $\bar{\lambda}(s)$, $\bar{\mu}(s)$ for the Darboux rectifying curves $\gamma(s)$ and $\bar{\gamma}(s)$ respectively does not change while taking an isometry on $M$ to $\bar{M}$.
\begin{theorem}
Let $f: M \rightarrow \bar{M}$ be an isometry of two smooth surfaces $M$ and $\bar{M}$. Then under the isometry, the component of the position vector of the Darboux rectifying curve along any tangent vector to the surface $M$ at $\gamma(s)$ is invariant, i.e.,$\gamma(s).T(=a\eta_{u} + b\eta_{v}) = \bar{\gamma}(s).\bar{T}(=a\bar{\eta_{u}} + b\bar{\eta_{v}})$.
\end{theorem}
\begin{proof}
Let $T=a\eta_{u} + b\eta_{v}$ be any tangent vector to the surface $M$ at $\gamma(s)$. Then
$$\gamma(s).T = [\lambda(s)(\eta_{u}u' + \eta_{v}v') + \frac{\mu(s)}{\sqrt{EG - F^{2}}}(Eu'\eta_{v} + F(v'\eta_{v} - u'\eta_{u}) - Gv'\eta_{u})].(a\eta_{u} + b\eta_{v}),$$
\begin{align*}
=\lambda(au'E + av'F + bu'F + bv'G) + \frac{\mu}{\sqrt{EG-F^{2}}}(EFau' + EGbu' + F^{2}au' \\+ FGbv' - FEau' - F^{2}bu' - GEav' - G^{2}bv'),
\end{align*}
\begin{align*}
= \lambda(au'E + F(av'+ bu') + bv'G) + \frac{\mu}{\sqrt{EG-F^{2}}}(F^{2}(av'- bu') + EG(bu'- av') \\+ FGbv' - G^{2}bv').
\end{align*}
On the other hand,
\begin{align*}
\bar{\gamma}(s).\bar{T}= \lambda(au'\bar{E} + \bar{F}(av'+ bu') + bv'\bar{G}) + \frac{\mu}{\sqrt{\bar{E}\bar{G}-\bar{F}^{2}}}(\bar{F}^{2}(av'- bu') \\+ \bar{E}\bar{G}(bu'- av') + \bar{F}\bar{G}bv' - \bar{G}^{2}bv').
\end{align*}
Using Theorem 3.1, where we get $\lambda(s) = \bar{\lambda}(s)$, $\mu(s) = \bar{\mu}(s)$ and under isometry, $E = \bar{E}$, $F = \bar{F}$ and $G = \bar{G}$. Hence $\gamma(s).T(=a\eta_{u} + b\eta_{v}) = \bar{\gamma}(s).\bar{T}(=a\bar{\eta_{u}} + b\bar{\eta_{v}})$.
\end{proof}
\begin{theorem} Let $f: M \rightarrow \bar{M}$ be an isometry between smooth surfaces $M$ and $\bar{M}$. Suppose $\gamma(s)$, $\bar{\gamma}(s)$ are Darboux rectifying curves on $M$, $\bar{M}$, respectively then for the component of $\gamma(s)$ and $\bar{\gamma}(s)$ along the principaal normals $N(s)$ and $\bar{N}(s)$ to the curves $\gamma(s)$ and $\bar{\gamma}(s)$ , the following relation holds:
$$\bar{\gamma}(s).\bar{N}(s) - \gamma(s).N(s) = \mu(s)A(E,F,G,E_{u},E_{v},F_{u},F_{v},G_{u},G_{v})(\frac{1}{\bar{\kappa}(s)}-\frac{1}{{\kappa}(s)}),$$
where $A(E,F,G,E_{u},E_{v},F_{u},F_{v},G_{u},G_{v}) = \frac{1}{\sqrt{EG-F^{2}}}[EG(u'v''-u''v') + F^{2}(u''v'-u'v'') + u'^{3}(E(F_{u}-\frac{E_{v}}{2})-F(\frac{E_{u}}{2})) + u'^{2}v'(EG_{u}-FE_{v}-G\frac{E_{u}}{2}+F(F_{u}-\frac{E_{v}}{2})) + u'v'^{2}(E(\frac{G_{v}}{2})+FG_{u}-F(F_{v}-\frac{G_{u}}{2})-GE_{v}) + v'^{3}(F(\frac{G_{v}}{2})-G(F_{v}-\frac{G_{u}}{2}))].$
\end{theorem}
\begin{proof}
Let $f: M \rightarrow \bar{M}$ be an isometry between smooth surfaces $M$ and $\bar{M}$ and $\gamma$ be a Darboux rectifying curve on $M$.
Then, we know that under the isometry between two smooth surfaces $M$ and $\bar{M}$,
\begin{equation}
E = \bar{E}, ~~~F = \bar{F}, ~~~G = \bar{G}.
\end{equation}
  Now from equation (3.3) \begin{equation}\bar{E_{u}} = E_{u}, \bar{E_{v}} = E_{v}, \bar{F_{u}} = F_{u}, \bar{F_{v}} = F_{v}, \bar{G_{u}} = G_{u}, \bar{G_{v}} = G_{v}.
 \end{equation} Also we know that
\begin{eqnarray}
\eta_{uu}.\eta_{u}=\frac{E_{u}}{2}, ~~~\eta_{uv}.\eta_{v}=\frac{G_{u}}{2}, ~~~\eta_{vv}.\eta_{v}=\frac{G_{v}}{2},\nonumber\\ \eta_{uv}.\eta_{u}=\frac{E_{v}}{2}, ~~~\eta_{uu}.\eta_{v}=F_{u}-\frac{E_{v}}{2}, ~~~\eta_{vv}.\eta_{u}=F_{v}-\frac{G_{u}}{2}.
\end{eqnarray}
Taking the component of $\gamma(s)$ along the normal to the curve, we get\\\\
$\gamma(s).N(s)= \frac{\mu(s)}{\kappa(s)\sqrt{EG - F^{2}}}(Eu'\eta_{v} + F(v'\eta_{v} - u'\eta_{u}) - Gv'\eta_{u}).(\eta_{u}u'' + \eta_{v}v'' + u'^{2}\eta_{uu} + 2u'v'\eta_{uv} + v'^{2}\eta_{vv})$.\\\\Then after simplification and using the equation (3.5), we obtain\\\\
$\gamma(s).N(s) = \frac{\mu(s)}{\kappa(s)\sqrt{EG - F^{2}}}[EG(u'v''-u''v') + F^{2}(u''v'-u'v'') + u'^{3}(E(F_{u}-\frac{E_{v}}{2})-F(\frac{E_{u}}{2})) + u'^{2}v'(EG_{u}-FE_{v}-G\frac{E_{u}}{2}+F(F_{u}-\frac{E_{v}}{2})) + u'v'^{2}(E(\frac{G_{v}}{2})+FG_{u}-F(F_{v}-\frac{G_{u}}{2})-GE_{v}) + v'^{3}(F(\frac{G_{v}}{2})-G(F_{v}-\frac{G_{u}}{2}))]$.\\\\
Similarly, we can also calculate $\bar{\gamma}(s).\bar{N}(s)$. Then using equations (3.3) and (3.4), we obtain\\
$$\bar{\gamma}(s).\bar{N}(s) - \gamma(s).N(s) = \mu(s)A(E,F,G,E_{u},E_{v},F_{u},F_{v},G_{u},G_{v})(\frac{1}{\bar{\kappa}(s)}-\frac{1}{{\kappa}(s)}),$$
where $A(E,F,G,E_{u},E_{v},F_{u},F_{v},G_{u},G_{v}) = \frac{1}{\sqrt{EG-F^{2}}}[EG(u'v''-u''v') + F^{2}(u''v'-u'v'') + u'^{3}(E(F_{u}-\frac{E_{v}}{2})-F(\frac{E_{u}}{2})) + u'^{2}v'(EG_{u}-FE_{v}-G\frac{E_{u}}{2}+F(F_{u}-\frac{E_{v}}{2})) + u'v'^{2}(E(\frac{G_{v}}{2})+FG_{u}-F(F_{v}-\frac{G_{u}}{2})-GE_{v}) + v'^{3}(F(\frac{G_{v}}{2})-G(F_{v}-\frac{G_{u}}{2}))].$
\end{proof}
\begin{corollary}
Let $f$ be an isometry of two smooth surfaces $M$ and $\bar{M}$. Then under the isometry, the component of the position vector of the Darboux rectifying curve along the normal vector to the curve at $\gamma(s)$ is invariant, i.e.,$\gamma(s).N(s) = \bar{\gamma}(s).\bar{N}(s)$ iff $\kappa(s) = \bar{\kappa}(s)$.
\end{corollary}
\begin{theorem}
Under the isometry $f: M \rightarrow \bar{M}$, the component of the position vector of the Darboux rectifying curve along binormal vector $B(s)$ to the curve at $\gamma(s)$, the following relation holds: $\bar{\gamma}(s).\bar{B}(s) - \gamma(s).B(s) = \mu(\frac{\bar{\kappa_{n}}(s)}{\bar{\kappa}(s)} - \frac{\kappa_{n}(s)}{\kappa(s)})$.
\end{theorem}
\begin{proof}
Let $f: M \rightarrow \bar{M}$ be an isometry between smooth surfaces $M$ and $\bar{M}$ and $\gamma$ be a Darboux rectifying curve on $M$. Then from equations (2.5) and (3.2), we obtain\\\\
$\gamma(s).B(s) = [\lambda(s)(\eta_{u}u' + \eta_{v}v') + \frac{\mu(s)}{\sqrt{EG - F^{2}}}(Eu'\eta_{v} + F(v'\eta_{v} - u'\eta_{u}) - Gv'\eta_{u})].[\frac{1}{\kappa(s)}\{(u'v''-v'u'')U + u'^{3}\eta_{u}\times \eta_{uu} + 2u'^{2}v'\eta_{u}\times \eta_{uv} + u'v'^{2}\eta_{u}\times \eta_{vv}  + u'^{2}v'\eta_{v}\times \eta_{uu} + 2u'v'^{2}\eta_{v}\times \eta_{uv} + v'^{3}\eta_{v}\times \eta_{vv}\}].$\\\\
After simplification, we get\\\\
$\gamma(s).B(s) = \frac{\mu}{\kappa\sqrt{EG-F^{2}}}[(Eu'^{4} + 2Fu'^{3}v' + Gu'^{2}v'^{2})\eta_{v}.(\eta_{u}\times \eta_{uu}) + (2Eu'^{3}v' + 4Fu'^{2}v'^{2} + 2Gu'v'^{3})\eta_{v}.(\eta_{u}\times \eta_{uv}) + (Eu'^{2}v'^{2} + 2Fu'v'^{3} + Gv'^{4})\eta_{v}.(\eta_{u}\times \eta_{vv})],$\\\\
$\Rightarrow\gamma(s).B(s) = \frac{\mu}{\kappa}[u'^{2}\eta_{uu}.U + 2u'v'\eta_{uv}.U + v'^{2}\eta_{vv}.U],$\\\\
$\Rightarrow\gamma(s).B(s) = \frac{\mu}{\kappa}[u'^{2}L + 2u'v'M + v'^{2}\mathbf{N}],$\\\\
$\Rightarrow\gamma(s).B(s) = \frac{\mu}{\kappa}\kappa_{n}$.\\\\
Similarly, we also get\\\\
$\bar{\gamma}(s).\bar{B}(s) = \frac{\mu}{\bar{\kappa}}\bar{\kappa_{n}}$.\\\\
Hence, $\bar{\gamma}(s).\bar{B}(s) - \gamma(s).B(s) = \mu(\frac{\bar{\kappa_{n}}(s)}{\bar{\kappa}(s)} - \frac{\kappa_{n}(s)}{\kappa(s)})$.
\end{proof}
\begin{corollary}
Let $f$ be an isometry of two smooth surfaces $M$ and $\bar{M}$. Then under the isometry, the component of the position vector of the Darboux rectifying curve along the binormal vector $B(s)$ to the curve at $\gamma(s)$ is invariant, i.e.,$\gamma(s).B(s) = \bar{\gamma}(s).\bar{B}(s)$ iff $\frac{\kappa_{n}}{\kappa} = \frac{\bar{\kappa_{n}}}{\bar{\kappa}}$.
\end{corollary}
\section{\textbf{Conformal image of a Darboux rectifying curve}}
\textbf{Definition 4.1.} {\em Let $\zeta$ be a diffeomorphism between two smooth surfaces $M$ and $\bar{M}$. Then $\zeta$ is said to be a local conformal map between $M$ and $\bar{M}$, if for all $y_{1},y_{2}\in T_{p}M$ with an arbitrary $p\in M$, $$<\zeta_{*}y_{1},\zeta_{*}y_{2}>_{\zeta(p)} = \rho^{2}<y_{1},y_{2}>_{p},$$ where $\rho$ is a differentiable function on M, also known as dilation factor.}\\\\
Therefore we can say that a conformal motion is the composition of an isometry and a dilation and when dilation factor is identity, then we get the isometry.
Now from \cite{mp},
\begin{equation}
\rho^{2}E = \bar{E}, ~~~~~\rho^{2}F = \bar{F}, ~~~~~\rho^{2}G = \bar{G}.\end{equation} That is, first fundamental form coefficients are conformally  invariant.\\\\
\textbf{Definition 4.2.} {\em Let $g: M \rightarrow \bar{M}$ be a smooth map between two smooth surfaces $M$ and $\bar{M}$ such that $\bar{g} = g\circ \eta$
with $\eta(u,v)$ being the surface patch of $M$. Then\\
i) $\bar{g}$ is said to conformally invariant when $\bar{g} = \rho^{2}g$ for some dilation factor $\rho(u,v)$.\\
ii) $\bar{g}$ is said to be homothetic invariant when $\bar{g} = c^{2}g$ for constant $c, c\neq \{0,1\}$.}
\begin{theorem}
Let $\zeta: M \rightarrow \bar{M}$ be a conformal map between two smooth surfaces $M$ and $\bar{M}$. Suppose $\gamma(s)$ be a Darboux rectifying curve on $M$ then $\bar{\gamma}(s)$ is a Darboux rectifying curve on $\bar{M}$ under $\zeta$ if $$\bar{\gamma}(s) = \zeta_{*}(\gamma(s)).$$
\end{theorem}
\begin{proof}
Let $\eta(u,v)$ and $\bar{\eta}(u,v) = \zeta\circ \eta(u,v)$ be the surface patches of $M$ and $\bar{M}$, respectively. The differential map $\zeta_{*}$ of $\zeta$ sends each tangent vector of $T_{p}M$ to a tangent vector of $T_{\zeta(p)}\bar{M}$. Then,
\begin{equation}
\bar{\eta_{u}}(u,v) = \zeta_{*}\eta_{u}, ~~~~~\bar{\eta_{v}}(u,v) = \zeta_{*}\eta_{v}.
\end{equation}
Suppose $\bar{\gamma}(s) = \zeta_{*}(\gamma(s))$. Then from the equation (3.2) of Darboux rectifying curve, we have
\begin{align*}
\Rightarrow \bar{\gamma}(s) = \lambda(s)(u'\zeta_{*}(\eta_{u}) + v'\zeta_{*}(\eta_{v})) + \frac{\mu(s)}{\sqrt{\rho^{2}E\rho^{2}G - \rho^{4}F^{2}}}\\(\rho^{2}Eu'\zeta_{*}\eta_{v} + \rho^{2}F(v'\zeta_{*}\eta_{v} - u'\zeta_{*}\eta_{u}) - \rho^{2}Gv'\zeta_{*}\eta_{u}).
\end{align*}
After using equations (4.1) and (4.2), we obtain
$$\bar{\gamma}(s) = \lambda(s)(u'\bar{\eta_{u}} + v'\bar{\eta_{v}}) + \frac{\mu(s)}{\sqrt{\bar{E}\bar{G} - \bar{F}^{2}}}(\bar{E}u'\bar{\eta_{v}} + \bar{F}(v'\bar{\eta_{v}} - u'\bar{\eta_{u}}) - \bar{G}v'\bar{\eta_{u}}),$$\\
$$\Rightarrow \bar{\gamma}(s) = \bar{\lambda}(s)\bar{T}(s) + \bar{\mu}(s)\bar{P}(s),$$\\
for smooth functions $\bar{\lambda}(s)$ and $\bar{\mu}(s)$. Thus $\bar{\gamma}(s)$ is a Darboux rectifying curve on $\bar{M}$.
\end{proof}
\textbf{Note:} From the Theorem 4.1., we see that the functions $\lambda(s)$, $\mu(s)$ and $\bar{\lambda}(s)$, $\bar\mu$(s) for the Darboux rectifying curves $\gamma(s)$ and $\bar{\gamma(s)}$, respectively does not change under the conformal map on $M$ to $\bar{M}$.
\begin{corollary}
Let $\zeta: M \rightarrow \bar{M}$ be a homothetic map and $\gamma(s)$ be a Darboux rectifying curve on $M$. Then $\bar{\gamma}(s)$ is a Darboux rectifying curve on $\bar{M}$ under $\zeta$ if $$\bar{\gamma}(s) = \zeta_{*}(\gamma(s)).$$
\end{corollary}
\begin{theorem}
Let $\zeta: M \rightarrow \bar{M}$ is a conformal map. Suppose $\gamma(s)$ and $\bar{\gamma}(s)$ are Darboux rectifying curve on $M$ and $\bar{M}$, respectively then for the component of the position vector of the Darboux rectifying curves along any tangent vector to the surface $M$ and $\bar{M}$ are given by $$\bar{\gamma}(s).\bar{T}(s) = \rho^{2}\gamma(s).T(s).$$
\end{theorem}
\begin{proof}
Let $T=a\eta_{u} + b\eta_{v}$ be any tangent vector to the surface $M$ at $\gamma(s)$. Then\\\\
$\rho^{2}\gamma(s).T(s) = \rho^{2}\{[\lambda(s)(\eta_{u}u' + \eta_{v}v') + \frac{\mu(s)}{\sqrt{EG - F^{2}}}(Eu'\eta_{v} + F(v'\eta_{v} - u'\eta_{u}) - Gv'\eta_{u})].(a\eta_{u} + b\eta_{v})\},$
\begin{align*}
\Rightarrow\rho^{2}\gamma(s).T(s) = \lambda(au'\rho^{2}E + \rho^{2}F(av'+ bu') + bv'\rho^{2}G) + \frac{\mu}{\sqrt{\rho^{2}E\rho^{2}G-\rho^{4}F^{2}}}\\(\rho^{4}F^{2}(av'- bu') + \rho^{2}E\rho^{2}G(bu'- av') + \rho^{2}F\rho^{2}Gbv' - \rho^{4}G^{2}bv').
\end{align*}
Then using equation (4.1), we obtain
\begin{align*}
\rho^{2}\gamma(s).T(s)= \lambda(au'\bar{E} + \bar{F}(av'+ bu') + bv'\bar{G}) + \frac{\mu}{\sqrt{\bar{E}\bar{G}-\bar{F}^{2}}}(\bar{F}^{2}(av'- bu') \\+ \bar{E}\bar{G}(bu'- av') + \bar{F}\bar{G}bv' - \bar{G}^{2}bv').
\end{align*}
Thus, $$\rho^{2}\gamma(s).T(s) = \bar{\gamma}(s).\bar{T}(s),$$\\
where we use the note of Theorem 4.1, $\lambda(s) = \bar{\lambda}(s)$ and $\mu(s) = \bar{\mu}(s)$ under the conformal map.
\end{proof}
\begin{theorem}
Let $\zeta: M \rightarrow \bar{M}$ be a conformal map and $\gamma(s)$ be a Darboux rectifying curve on $M$, then for the component of $\gamma(s)$ along the principal normal $N(s)$ to the curve $\gamma(s)$ , the following relation holds:
$$\rho^{2}\gamma(s).N(s) - \bar{\gamma}(s).\bar{N}(s) =  \psi(E,F,G,\rho)$$ iff $\kappa(s) = \bar{\kappa}(s)$,
where $\psi(E,F,G,\rho) = \frac{\mu(s)}{\kappa(s)\sqrt{EG - F^{2}}}[u'^{3}(E^{2}\rho\rho_{v} - EF\rho\rho_{u}) + u'^{2}v'(3EF\rho\rho_{v} - EG\rho\rho_{u} - 2F^{2}\rho\rho_{u}) + u'v'^{2}(2F^{2}\rho\rho_{v} - 3FG\rho\rho_{u} + EG\rho\rho_{v}) + v'^{3}(FG\rho\rho_{v} - G^{2}\rho\rho_{u})].$
\end{theorem}
\begin{proof}
Let $\zeta: M \rightarrow \bar{M}$ be a conformal map and $\gamma(s)$ be a Darboux rectifying curve on $M$. Then from equation (4.1), we obtain\\
\begin{eqnarray}\bar{E_{u}} = 2\rho\rho_{u}E + \rho^{2}E_{u}, \bar{E_{v}} = 2\rho\rho_{v}E + \rho^{2}E_{v},\nonumber\\ \bar{F_{u}} = 2\rho\rho_{u}F + \rho^{2}F_{u}, \bar{F_{v}} = 2\rho\rho_{v}F + \rho^{2}F_{v},\nonumber\\ \bar{G_{u}} = 2\rho\rho_{u}G + \rho^{2}G_{u}, \bar{G_{v}} = 2\rho\rho_{v}G + \rho^{2} G_{v}.
\end{eqnarray}
 Also we know that
\begin{eqnarray}
\eta_{uu}.\eta_{u}=\frac{E_{u}}{2}, ~~~\eta_{uv}.\eta_{v}=\frac{G_{u}}{2}, ~~~\eta_{vv}.\eta_{v}=\frac{G_{v}}{2},\nonumber\\ \eta_{uv}.\eta_{u}=\frac{E_{v}}{2}, ~~~\eta_{uu}.\eta_{v}=F_{u}-\frac{E_{v}}{2}, ~~~\eta_{vv}.\eta_{u}=F_{v}-\frac{G_{u}}{2}.
\end{eqnarray}
Taking the component of $\gamma(s)$ along the normal to the curve, we get\\\\
$\gamma(s).N(s)= \frac{\mu(s)}{\kappa(s)\sqrt{EG - F^{2}}}(Eu'\eta_{v} + F(v'\eta_{v} - u'\eta_{u}) - Gv'\eta_{u}).(\eta_{u}u'' + \eta_{v}v'' + u'^{2}\eta_{uu} + 2u'v'\eta_{uv} + v'^{2}\eta_{vv})$.\\\\Then after simplification and using the equation (4.3) and (4.4), we obtain\\\\
$\rho^{2}\gamma(s).N(s) = \frac{\rho^{4}\mu(s)}{\kappa(s)\sqrt{\rho^{2}E\rho^{2}G - \rho^{4}F^{2}}}[EG(u'v''-u''v') + F^{2}(u''v'-u'v'') + u'^{3}(E(F_{u}-\frac{E_{v}}{2})-F(\frac{E_{u}}{2})) + u'^{2}v'(EG_{u}-FE_{v}-G\frac{E_{u}}{2}+F(F_{u}-\frac{E_{v}}{2})) + u'v'^{2}(E(\frac{G_{v}}{2})+FG_{u}-F(F_{v}-\frac{G_{u}}{2})-GE_{v}) + v'^{3}(F(\frac{G_{v}}{2})-G(F_{v}-\frac{G_{u}}{2}))]$.\\\\
$\Rightarrow \rho^{2}\gamma(s).N(s) = \frac{\mu(s)}{\kappa(s)\sqrt{\rho^{2}E\rho^{2}G - \rho^{4}F^{2}}}[\rho^{2}E\rho^{2}G(u'v''-u''v') + \rho^{4}F^{2}(u''v'-u'v'') + u'^{3}(\rho^{2}E(\bar{F_{u}}-\frac{\bar{E_{v}}}{2})-\rho^{2}F(\frac{\bar{E_{u}}}{2})) + u'^{2}v'(\rho^{2}E\bar{G_{u}}-\rho^{2}F\bar{E_{v}}-\rho^{2}G\frac{\bar{E_{u}}}{2}+\rho^{2}F(\bar{F_{u}}-\frac{E_{v}}{2})) + u'v'^{2}(\rho^{2}E(\frac{\bar{G_{v}}}{2})+\rho^{2}F\bar{G_{u}}-\rho^{2}F(\bar{F_{v}}-\frac{\bar{G_{u}}}{2})-\rho^{2}G\bar{E_{v}}) + v'^{3}(\rho^{2}F(\frac{\bar{G_{v}}}{2})-\rho^{2}G(\bar{F_{v}}-\frac{\bar{G_{u}}}{2}))] + \frac{\mu(s)}{\kappa(s)\sqrt{EG - F^{2}}}[u'^{3}(E^{2}\rho\rho_{v} - EF\rho\rho_{u}) + u'^{2}v'(3EF\rho\rho_{v} - EG\rho\rho_{u} - 2F^{2}\rho\rho_{u}) + u'v'^{2}(2F^{2}\rho\rho_{v} - 3FG\rho\rho_{u} + EG\rho\rho_{v}) + v'^{3}(FG\rho\rho_{v} - G^{2}\rho\rho_{u})],$\\\\
$\Rightarrow \rho^{2}\gamma(s).N(s) = \frac{\mu(s)}{\kappa(s)\sqrt{\bar{E}\bar{G} - \bar{F}^{2}}}[\bar{E}\bar{G}(u'v''-u''v') + \bar{F}^{2}(u''v'-u'v'') + u'^{3}(\bar{E}(\bar{F_{u}}-\frac{\bar{E_{v}}}{2})-\bar{F}(\frac{\bar{E_{u}}}{2})) + u'^{2}v'(\bar{E}\bar{G_{u}}-\bar{F}\bar{E_{v}}-\bar{G}\frac{\bar{E_{u}}}{2}+\bar{F}(\bar{F_{u}}-\frac{E_{v}}{2})) + u'v'^{2}(\bar{E}(\frac{\bar{G_{v}}}{2})+\bar{F}\bar{G_{u}}-\bar{F}(\bar{F_{v}}-\frac{\bar{G_{u}}}{2})-\bar{G}\bar{E_{v}}) + v'^{3}(\bar{F}(\frac{\bar{G_{v}}}{2})-\bar{G}(\bar{F_{v}}-\frac{\bar{G_{u}}}{2}))] + \psi(E,F,G,\rho),$\\\\
$\Rightarrow \rho^{2}\gamma(s).N(s) = \bar{\gamma}(s).\bar{N}(s) +  \psi(E,F,G,\rho),$\\\\
$\Rightarrow \rho^{2}\gamma(s).N(s) - \bar{\gamma}(s).\bar{N}(s) =  \psi(E,F,G,\rho),$\\\\
where $\psi(E,F,G,\rho) = \frac{\mu(s)}{\kappa(s)\sqrt{EG - F^{2}}}[u'^{3}(E^{2}\rho\rho_{v} - EF\rho\rho_{u}) + u'^{2}v'(3EF\rho\rho_{v} - EG\rho\rho_{u} - 2F^{2}\rho\rho_{u}) + u'v'^{2}(2F^{2}\rho\rho_{v} - 3FG\rho\rho_{u} + EG\rho\rho_{v}) + v'^{3}(FG\rho\rho_{v} - G^{2}\rho\rho_{u})].$ Hence the proof.
\end{proof}
\begin{theorem}
Let $\zeta$ be a conformal map between two smooth surfaces $M$ and $\bar{M}$. Suppose $\gamma(s)$ be a Darboux rectifying curve on $M$ then for the component of $\gamma(s)$ along the binormal $B(s)$ to the curve at $\gamma(s)$, the following relation holds:\\\\
$ \bar{\gamma}(s).\bar{B}(s) - \rho^{2}\gamma(s).B(s) =  \frac{\mu}{\kappa}[u'^{2}(\frac{\bar{\eta_{uu}}}{\bar{W}^{2}} - \rho^{2}\frac{\eta_{uu}}{W^{2}}) + 2u'v'(\frac{\bar{\eta_{uv}}}{\bar{W}^{2}} - \rho^{2}\frac{\eta_{uv}}{W^{2}}) + v'^{2}(\frac{\bar{\eta_{vv}}}{\bar{W}^{2}} - \rho^{2}\frac{\eta_{vv}}{W^{2}})],$ where $W^{2} = 1 + \eta_{u}^{2} + \eta_{v}^{2},$ $\bar{W}^{2} = 1 + \bar{\eta_{u}}^{2} + \bar{\eta_{v}}^{2}$, iff $\bar{\kappa} =  \kappa$
\end{theorem}
\begin{proof}
Let $\zeta$ be a conformal map between two smooth surfaces $M$ and $\bar{M}$ and $\gamma(s)$ be a Darboux rectifying curve on $M$. Then from the equation (2.5) and (3.2), we obtain\\\\
$\gamma(s).B(s) = [\lambda(s)(\eta_{u}u' + \eta_{v}v') + \frac{\mu(s)}{\sqrt{EG - F^{2}}}(Eu'\eta_{v} + F(v'\eta_{v} - u'\eta_{u}) - Gv'\eta_{u})].[\frac{1}{\kappa(s)}\{(u'v''-v'u'')U + u'^{3}\eta_{u}\times \eta_{uu} + 2u'^{2}v'\eta_{u}\times \eta_{uv} + u'v'^{2}\eta_{u}\times \eta_{vv}  + u'^{2}v'\eta_{v}\times \eta_{uu} + 2u'v'^{2}\eta_{v}\times \eta_{uv} + v'^{3}\eta_{v}\times \eta_{vv}\}].$\\\\
After simplification, we obtain\\\\
$\gamma(s).B(s) = \frac{\mu}{\kappa}[u'^{2}L + 2u'v'M + v'^{2}\mathbf{N}]$.\\\\
Thus\\\\
$\bar{\gamma}(s).\bar{B}(s) - \rho^{2}\gamma(s).B(s) = \mu[u'^{2}(\frac{\bar{L}}{\bar{\kappa}} - \rho^{2}\frac{L}{\kappa}) + 2u'v'(\frac{\bar{M}}{\bar{\kappa}} - \rho^{2}\frac{M}{\kappa}) + v'^{2}(\frac{\bar{\mathbf{N}}}{\bar{\kappa}} - \rho^{2}\frac{\mathbf{N}}{\kappa})].$\\\\
Now in the Monge patch form, $L$, $M$ and $\mathbf{N}$, can written in the given form\\\\
$L = \frac{\eta_{uu}}{1 + \eta_{u}^{2} + \eta_{v}^{2} (= W^{2})}$, $M = \frac{\eta_{uv}}{1 + \eta_{u}^{2} + \eta_{v}^{2} (= W^{2})}$ and $\mathbf{N} = \frac{\eta_{vv}}{1 + \eta_{u}^{2} + \eta_{v}^{2} (= W^{2})}$, and similar expression for $\bar{L}, \bar{M}, \bar{\mathbf{N}}$. Hence\\\\
$ \bar{\gamma}(s).\bar{B}(s) - \rho^{2}\gamma(s).B(s) =  \frac{\mu}{\kappa}[u'^{2}(\frac{\bar{\eta_{uu}}}{\bar{W}^{2}} - \rho^{2}\frac{\eta_{uu}}{W^{2}}) + 2u'v'(\frac{\bar{\eta_{uv}}}{\bar{W}^{2}} - \rho^{2}\frac{\eta_{uv}}{W^{2}}) + v'^{2}(\frac{\bar{\eta_{vv}}}{\bar{W}^{2}} - \rho^{2}\frac{\eta_{vv}}{W^{2}})],$ where $W^{2} = 1 + \eta_{u}^{2} + \eta_{v}^{2},$ $\bar{W}^{2} = 1 + \bar{\eta_{u}}^{2} + \bar{\eta_{v}}^{2}$, iff $\bar{\kappa} =  \kappa$

\end{proof}

\end{document}